\documentclass[a4paper,10pt]{article}

\usepackage{amsmath,amsfonts,amssymb,latexsym,amsthm, enumerate, verbatim, a4wide}

  \theoremstyle{plain}
  \newtheorem{thm}{Theorem}[section]
  \newtheorem{cor}[thm]{Corollary}
  \newtheorem{prop}[thm]{Proposition}
  \newtheorem{lem}[thm]{Lemma}
  \newtheorem{defi}[thm]{Definition}

\DeclareMathOperator{\ric}{Ric}

\DeclareMathOperator{\vol}{vol}

\DeclareMathOperator{\bino}{bin}

\DeclareMathOperator{\supp}{Supp}

\DeclareMathOperator{\EG}{E\Gamma}

\DeclareMathOperator{\diver}{div}
\DeclareMathOperator{\grad}{grad}
\title{Entropy along $W_{1,+}$-geodesics on graphs}

\author{Erwan Hillion \footnote{Department of Mathematics, University of Luxembourg, erwan.hillion@uni.lu}}

\begin{document}

\maketitle 

\begin{abstract}
We study the convexity of the entropy functional along particular interpolating curves defined on the space of finitely supported probability measures on a graph. 
\end{abstract}

\section{Introduction}

The Wasserstein distance $W_p(\mu_0,\mu_1)$ between two finitely supported probability measures on a metric space $(X,d)$ with its Borel $\sigma$-algebra is defined for $p \geq 1$ by
\begin{equation} \label{eq:WassersteinDef}
W_p(\mu_0,\mu_1)^p := \inf_{\pi \in \Pi(\mu_0,\mu_1)} \int_{X \times X} d(x_0,x_1)^p d\pi(x_0,x_1),
\end{equation}
where $\Pi(\mu_0,\mu_1)$ is the (non-empty) set of couplings between $\mu_0$ and $\mu_1$, i.e. the set of probability measures on $X \times X$ having $\mu_0$ and $\mu_1$ as marginals. The optimization problem defined by equation~\eqref{eq:WassersteinDef} is called the Monge-Kantorovitch problem and any minimizer for~\eqref{eq:WassersteinDef} is called optimal coupling between $\mu_0$ and $\mu_1$. For a comprehensive study of optimal transportation theory, the reader is referred to the textbooks~\cite{VillaniBook1} and~\cite{VillaniBook2} by Villani.

\medskip

Under mild conditions, it is possible to show that the set $\Pi_p(\mu_0,\mu_1)$ of optimal couplings between $\mu_0$ and $\mu_1$ is non-empty. Furthermore, under the additional assumptions that $p>1$, $(X,d)$ is the Euclidean space $(\mathbb{R}^d,|.|)$ and $\mu_0$ is absolutely continuous with respect to the Lebesgue measure, one can prove the existence of a measurable map $T : \mathbb{R}^d \rightarrow \mathbb{R}^d$ such that the coupling $\pi := (Id \times T)_*\mu_0$ is a minimizer for~\eqref{eq:WassersteinDef}.
 
\medskip 
 
In particular, $\mu_1$ is the pushforward of $\mu_0$ by the application $T$: $\mu_1 := T_*\mu_0$ and equation~\eqref{eq:WassersteinDef} can be rewritten
\begin{equation}
W_p(\mu_0,\mu_1)^p = \int_{\mathbb{R}^d} |x-T(x)|^p d\mu_0(x).
\end{equation}
It is possible to go further by considering, for $0 \leq t \leq 1$, the measure $\mu_t := (T_t)_*\mu_0$, where the application $T_t : \mathbb{R}^d \rightarrow \mathbb{R}^d$ is defined as the barycenter $T_t(x) := (1-t)x + t T(x)$. One can then show that the family $(\mu_t)_{t \in [0,1]}$ is a geodesic for the Wasserstein distance $W_p$, in the sense that 
\begin{equation*}
W_p(\mu_0,\mu_1) = \sup_{0 = t_0 \leq t_1 \cdots \leq t_{n} = 1} \sum_{i=0}^{n-1} W_p(\mu_{t_i},\mu_{t_{i+1}}).
\end{equation*}
Moreover, a fundamental property of optimal couplings asserts that $T_t$ is injective, which allows us to define unambiguously a velocity field $(v_t)_{t \in [0,1]}$ by 
\begin{equation*}
v_t(T_t(x)) := T(x)-x.
\end{equation*}

The terminology 'velocity field' comes from the fact that, if we write $d\mu_t(x) = f_t(x) dx$, then the density $f_t(x)$ satisfy, at least formally, the transport equation
\begin{equation} \label{eq:velocitytransport}
\frac{\partial}{\partial t} f_t(x) + \diver (v_t(x) f_t(x)) = 0. 
\end{equation}

Moreover, the velocity field $v_t(x)$ satisfies the Hamilton-Jacobi-type equation
\begin{equation} \label{eq:velocityoptimalHJ}
\frac{\partial}{\partial t} v_t(x) + \frac{1}{2} \grad |v_t(x)|^2 = 0,
\end{equation}
which can be simplified into
\begin{equation} \label{eq:velocityoptimal}
\frac{\partial}{\partial t} v_t(x) = - \diver (v_t(x)) v_t(x). 
\end{equation}

In~\cite{BenamouBrenier}, Benamou and Brenier proved that both equations~\eqref{eq:velocitytransport} and~\eqref{eq:velocityoptimal} can be used to give a characterization of $W_p$-geodesics, more precisely we have:
\begin{thm}\label{th:BenamouBrenier}
Given two finitely supported probability measures $d\mu_0(x) := f_0(x) dx$ and $d\mu_1(x) := f_1(x) dx$, we have
\begin{equation} \label{eq:MKBenamou}
W_p(\mu_0,\mu_1)^p = \inf \int_0^1 \int_{\mathbb{R}^d} |v_t(x)|^p d\mu_t(x),
\end{equation}
where the infimum is taken over the set of curves $(\mu_t)_{t \in [0,1]} = (f_t(x) dx)_{t \in [0,1]}$ joining the prescribed measures $\mu_0$ and $\mu_1$, and where $(v_t(x))_{t \in [0,1]}$ is a velocity field such that equation~\eqref{eq:velocitytransport} holds. Moreover, the formal optimality condition for the optimization problem~\eqref{eq:MKBenamou} is given by equation~\eqref{eq:velocityoptimalHJ}.
\end{thm}

Theorem~\ref{th:BenamouBrenier} is also true for families of probability measures defined on a Riemannian manifold, having smooth enough densities with respect to the Riemannian volume measure. However, in this framework, equations \eqref{eq:velocityoptimalHJ} and \eqref{eq:velocityoptimal} are no longer equivalent.

\medskip

The optimality condition~\eqref{eq:velocityoptimal} is the starting point of the article~\cite{HillionGeodesic} by the author. The main idea is the following: given two distinct probability measures $f_0,f_1$ on a graph $G$, there is no interpolating curve $(f_t)_{t \in [0,1]}$ with a finite length for the Wasserstein $W_p$, for any $p>1$. However, in generic cases there are infinitely many geodesics $(f_t)_{t \in [0,1]}$ for the $W_1$ distance. The aim of~\cite{HillionGeodesic} is to choose among this set a particular $W_1$-geodesic satisfying a discrete version of equation~\eqref{eq:velocityoptimal}. These interpolating curves are called $W_{1,+}$-geodesics on $G$; we recall their basic properties in Section 2.

\medskip

The purpose of this article is to study the behaviour of the entropy functional along a $W_{1,+}$-geodesic $(f_t(x))_{t \in [0,1], x \in G}$ on a graph $G$. More precisely, we will study the convexity of the function $t \mapsto H(t)$ defined by
\begin{equation}
H(t) := \sum_{x \in G} f_t(x) \log(f_t(x)), 
\end{equation}
where by convention $0 \log 0 = 0$. The methods used to prove such convexity properties are adapted from the previous article~\cite{HillionContraction} by the author, and use the first-order-calculus formalism introduced in~\cite{HillionGeodesic}.

\medskip

The motivation behind this research work comes from Sturm-Lott-Villani theory, developed in the articles~\cite{SturmRicci01}, \cite{SturmRicci02} and ~\cite{LottVillani}. The main idea of this theory is the following: it is possible to obtain some information about the geometry of a measured length space $(X,d,\nu)$ by studying the behaviour of entropy functionals along $W_2$-geodesics on the space of probability measures over $(X,d)$. A major result asserts that a compact Riemannian manifold $(M,g)$ satisfies the Ricci curvature bound $\ric \geq K g$ if and only if each pair of absolutely continuous probability measures $\mu_0, \mu_1$ can be joined by a Wasserstein $W_2$-geodesic $(\mu_t)_{t \in [0,1]}$ such that 
\begin{equation}\label{eq:LVSconvexity}
H(\mu_t) \leq (1-t) H(\mu_0) + t H(\mu_1) - K \frac{t(1-t)}{2} W_2(\mu_0,\mu_1)^2,
\end{equation}
where the relative entropy $H(\mu)$ is defined by $H(\mu): = \int_{M} \rho \log(\rho) d\vol$ if $d\mu = \rho . d\vol$ and by $H(\mu)=\infty$ if $\mu$ is not absolutely continuous with respect to the Riemannian volume measure. It is then possible to define the curvature condition '$\ric \geq K$' on a measured length space $(X,d,\nu)$ if Equation~\eqref{eq:LVSconvexity} is satisfied for any $W_2$-Wasserstein geodesic on $\mathcal{P}_2(X)$. Several geometric theorems and functional inequalities holding on Riemannian manifolds satisfying a Ricci curvature bound are still valid in the framework of measured length spaces with a curvature condition '$\ric \geq K$'.

\medskip

The generalization of Sturm-Lott-Villani theory to discrete setting has been the subject of many research works, each leading to its own definition of Ricci curvature bounds on graphs, among which we can cite papers by Ollivier~\cite{OllivierRicci} and Erbar-Maas \cite{ErbarMaas}. The latter is based on the study of a discrete version of the minimization problem~\eqref{eq:MKBenamou} for $p=2$, whereas our approach is based on a discrete version of equation~\eqref{eq:velocityoptimal} characterizing the solutions of~\eqref{eq:MKBenamou}. Another important work in discrete Sturm-Lott-Villani theory is~\cite{GRST} which, like this present work, is based on the study of the behaviour of the entropy functional along mixtures of binomial measures. 

\medskip

The results proven in our paper show that the convexity properties of the entropy along $W_{1,+}$-geodesics are linked with some intuitive notion of curvature bounds on graphs. However, it seems that our study of the convexity of the entropy does not lead to a  definition of Ricci curvature bounds strong enough to imply important functional inequalities, such as the modified logarithmic Sobolev inequality introduced in~\cite{BobkovLedoux}. 

\medskip

Our article is outlined as follows: in Section~\ref{sec:Rappel}, we recall the definition and basic properties of the $W_1$-orientation and $W_{1,+}$-geodesics, which are developed in the previous article~\cite{HillionGeodesic}. We also introduce the notion of canonical $W_{1,+}$-geodesic, see Theorem~\ref{th:CanonicalGeodesic}, which will be used in Section~\ref{sec:ProductSpaces}. 

\medskip

In Section~\ref{sec:GeneralFormula}, we begin the study of the entropy function $H(t)$ along a $W_{1,+}$-geodesic on a graph; we use the Benamou-Brenier equation~\eqref{eq:BBcondition}, which is at the heart of definition of $W_{1,+}$-geodesics, to obtain bounds on the second derivative $H''(t)$. The calculations done in this section are inspired by those done in the previous article~\cite{HillionContraction} by the author; the results obtained are also linked with the more general theory of entropic interpolations, developed by L\'eonard in a recent series of articles, including~\cite{LeoConvex}, \cite{LeoLazy} and \cite{LeonardSurvey}. 

\medskip

In Section~\ref{sec:ProductSpaces}, we refine the calculations done in Section~\ref{sec:GeneralFormula} to prove a 
tensorization property. This property allows us to give bounds on the second derivative $H''(t)$ when the underlying graph is a product graph. Interesting examples are given by $\mathbb{Z}^n$, the cube $\{0,1\}^n$, or more generally by the Cayley graph of a finitely generated abelian group. 

\medskip

In the Appendix, we present two additional results on families of probability measures on $\mathbb{Z}$. We first prove that, along a $W_{1,+}$-geodesic on $\mathbb{Z}$, other types of functionals are convex, belonging to the family of Renyi entropies functionals. The second part of the Appendix is devoted to another type of interpolation of probability measures on $\mathbb{Z}$, defined as a mixture of binomial distributions with respect to a $W_2$-optimal coupling.

\section{$W_{1,+}$-geodesics on graphs}\label{sec:Rappel}

In this section, we first recall the main definitions and properties of \cite{HillionGeodesic}. The reader is referred to this paper for detailed proofs and additional explanations. We then introduce the new notion of canonical $W_{1,+}$-geodesic, which will be used in the study of product spaces in Section~\ref{sec:ProductSpaces}.

\subsection{Definition and construction}

Let $G$ be a locally finite, connected graph. We denote by $d$ the usual graph distance on $G$ and by $x \sim y$ the adjacency relation on $G$, meaning that $(x,y)$ is an edge of $G$. A curve of length $n$ on $G$ is an application $\gamma : \{0,\ldots n\} \rightarrow G$ satisfying $\gamma(i) \sim \gamma(i+1)$. A geodesic between two vertices $x$ and $y$ is a curve of minimal length joining $x$ to $y$. The set of geodesics between $x$ and $y$ is denoted by $\Gamma_{x,y}$ and its cardinality by $|\Gamma_{x,y}|$. The set of all geodesic curves of $G$ is denoted by $\Gamma(G)$.

\medskip

Let $f_0,f_1$ be two finitely supported probability distributions on $G$. We denote by $\Pi_1(f_0,f_1)$ the set of $W_1$-optimal couplings between $f_0$ and $f_1$, i.e. the set of couplings between $f_0$ and $f_1$ which minimize the functional
\begin{equation*}
I_1(\pi) := \sum_{x,y \in G} d(x,y) \pi(x,y).
\end{equation*}

Using properties of supports of optimal couplings, one can prove that the following definition in unambiguous:
\begin{defi} Let $f_0,f_1$ be two finitely supported probability measures on $G$.
\begin{itemize}
\item The $W_1$-orientation on $G$ with respect to $f_0,f_1$ is constructed in the following way: a couple $(x,y)$ of adjacent vertices is oriented by $x \rightarrow y$ if there exists an optimal coupling $\pi \in \Pi_1(f_0,f_1)$ and a geodesic $\gamma \in \Gamma(G)$ of length $n$ such that $(\gamma(0),\gamma(n)) \in \supp(\pi)$ and such that there exists $i \in \{0,\ldots n-1\}$ with $\gamma(i)=x$ and $\gamma(i+1)=y$.
\item Let $x_1 \in G$. The set $\mathcal{E}(x_1)$, resp. $\mathcal{F}(x_1)$, is the (possibly empty) set of vertices $x_0 \in G$, resp. $x_2 \in G$, such that $x_0 \rightarrow x_1$, resp. $x_1 \rightarrow x_2$.
\item An oriented path on $G$ is a mapping $\gamma : \{0,\ldots n\}$ with $\gamma(i) \rightarrow \gamma(i+1)$.
\item The $W_1$-orientation w.r.t. $f_0,f_1$ induces a partial order on the vertices of $G$: we denote $x \leq y$ if there exists an oriented path $x = \gamma(0) \rightarrow \cdots \rightarrow \gamma(n)=y$.
\end{itemize}
\end{defi}

One important property of this orientation is the fact that every oriented path is a geodesic:

\begin{prop}\label{prop:OrientedGeodesic}
If we have $\gamma(0) \rightarrow \cdots \rightarrow \gamma(n)$ then $d(\gamma(0),\gamma(n)) = n$.
\end{prop}

A particular subset of geodesics on the oriented $G$ is given by extremal geodesics:

\begin{defi}
Let $\gamma : \gamma(0) \rightarrow \cdots \rightarrow \gamma(n)$ be a geodesic on the oriented $G$. We say that $\gamma$ is an extremal geodesic, and we write $\gamma \in \EG$ if it cannot be extended in a longer geodesic, i.e. if the sets $\mathcal{E}(\gamma(0))$ and $\mathcal{F}(\gamma(n))$ are empty.
\end{defi}
\medskip

The introduction of an orientation makes possible the introduction of a first-order calculus on $G$. We first define:

\begin{defi}
The oriented edge graph $(E(G),\rightarrow)$ associated to $(G,\rightarrow)$ is defined as follows: its vertices are denoted by $(x_0x_1)$, where $x_0 \rightarrow x_1 \in G$ and its oriented edges join each couple $(x_0x_1) \rightarrow (y_0y_1)$ such that $x_1=y_0$.

\medskip

The oriented graph of oriented triples $(T(G),\rightarrow)$ is the graph $(E(E(G)),\rightarrow)$: its vertices are the triples $(x_0x_1x_2)$ with $x_0 \rightarrow x_1 \rightarrow x_2$ and its edges are defined between each couple $(x_0x_1x_2) \rightarrow (x_1x_2x_3)$.
\end{defi}

When the choice of the orientation on $G$ is unambiguous, we will often write $E(G), T(G)$ instead of $(E(G),\rightarrow), (T(G),\rightarrow)$.

\begin{defi}
The divergence of a function $g : E(G) \rightarrow \mathbb{R}$ defined on the oriented edges of $G$ is the function $\nabla \cdot  g : G \rightarrow  \mathbb{R}$ defined by:
\begin{equation}
\nabla \cdot  g (x_1)= \sum_{x_2 \in \mathcal{F}(x_1)} g(x_1x_2) - \sum_{x_0 \in \mathcal{E}(x_1)} g(x_0x_1).
\end{equation}
We define similarly the divergence $\nabla \cdot  h : E(G) \rightarrow \mathbb{R}$ of a function $h : T(G) \rightarrow \mathbb{R}$ defined on the oriented triples of $G$. We denote $\nabla_2 \cdot h := \nabla \cdot (\nabla \cdot h)$.
\end{defi}

\medskip

This first-order differential operator on the oriented graph allows us to introduce a discrete version of the formal optimality condition~\eqref{eq:velocityoptimal}, on which is based the definition of $W_{1,+}$-geodesics:

\begin{defi}
Let $G$ be a graph, $W_1$-oriented with respect to a couple of probability measures $f_0,f_1$. A family $(f_t)=(f_t)_{t \in [0,1]}$ is said to be a $W_{1,+}$-geodesic if:
\begin{enumerate}
\item The curve $(f_t)$ is a $W_1$-geodesic.
\item There exist two families $(g_t)$ and $(h_t)$ defined respectively on $E(G)$ and $T(G)$, such that:
\begin{equation*}
\frac{\partial}{\partial t} f_t = - \nabla \cdot  g_t \ , \ \frac{\partial}{\partial t} g_t = - \nabla \cdot  h_t.
\end{equation*}
\item For every $(xy) \in E(G)$ we have $g_t(xy)>0$.
\item The triple $(f_t,g_t,h_t)$ satisfies the Benamou-Brenier equation
\begin{equation}\label{eq:BBcondition}
\forall (x_0x_1x_2) \in T(G) \ , \ f_t(x_1)h_t(x_0x_1x_2) = g_t(x_0x_1) g_t(x_1x_2).
\end{equation}
\end{enumerate}
\end{defi}

\medskip

Let us fix a couple $f_0$, $f_1$ of probability measures on $G$ and endow $G$ with the $W_1$-orientation with respect to $f_0,f_1$. The existence of a $W_{1,+}$-interpolation $(f_t)$ joining $f_0$ to $f_1$ is the main result of~\cite{HillionGeodesic}. Moreover, any such curve $(f_t)$ can be seen as a mixture of binomial families of distributions with respect to a coupling which is solution of a certain minimization problem. 

\subsection{Canonical $W_{1,+}$-geodesics}

In this paper we are mostly interested in particular $W_{1,+}$-geodesics, called canonical $W_{1,+}$-geodesics, which correspond to the case where $\forall  \gamma \in \EG \ , \ C(\gamma)=1$, with the notations of~\cite{HillionGeodesic}. The existence, uniqueness, and construction of such curves can be summed up by the following:

\begin{thm} \label{th:CanonicalGeodesic}
Let $x_0 \leq \cdots \leq x_n \in G$ be an oriented $n+1$-uples of vertices of $G$. We define: 
\begin{equation*} 
m(x_0,\ldots, x_n) := \frac{|\Gamma_{x_0,\cdots, x_n}|}{|\EG|},
\end{equation*}
where $\Gamma_{x_0,\cdots, x_n}$ is the set of extremal geodesics visiting $x_0,\ldots, x_n$: $$\Gamma_{x_0,\cdots x_n} := \{\gamma \in \EG \ :  \ \exists i_0 \leq \cdots \leq i_n \ ,\ \gamma(i_k)= x_k\}.$$
There exists a unique couple of families of functions $P_t(x),Q_t(x)$, defined for $x \in G$ and $t \in [0,1]$, such that each $t \mapsto P_t(x)$ and $t \mapsto Q_t(x)$ is positive and polynomial in $t$, and satisfying the following property: let us consider the families of functions $(f_t)$, $(g_t)$, $(h_t)$ respectively defined on $G$, $E(G)$ and $T(G)$ by \begin{eqnarray*} f_t(x_0) &:=& m(x_0)P_t(x_0)Q_t(x_0), \\  g_t(x_0x_1) &:=& m(x_0,x_1) P_t(x_0) Q_t(x_1), \\ h_t(x_0x_1x_2) &:=& m(x_0,x_1,x_2) P_t(x_0)Q_t(x_2). \end{eqnarray*} Then the triple $(f_t,g_t,h_t)$ satisfies all the items of the definition of a $W_{1,+}$-geodesic. Such a curve will be called canonical $W_{1,+}$-geodesic joining $f_0$ to $f_1$.
\end{thm}

The reason why we introduce these particular geodesics comes from the following property, which will be used in Section~\ref{sec:ProductSpaces}:

\begin{prop} \label{prop:CanonicalTriple}
If the triple $(f_t,g_t,h_t)$ defines a canonical $W_{1,+}$-geodesic, then for any oriented triple $(x_0x_1x_2) \in T(G)$, the quantity $h(x_0x_1x_2)$ does not depend on $x_1$, and therefore can be written $h(x_0x_2)$.
\end{prop}

\textbf{Proof:} It suffices to show that the cardinality $|\Gamma_{x_0,x_1,x_2}|$ does not depend on $x_1$. This comes from the fact that every $\gamma \in \Gamma_{x_0,x_1,x_2}$ can be written $$ \gamma \ : \  \gamma(0) \rightarrow \cdots \rightarrow \gamma(i)=x_0 \rightarrow x_1 \rightarrow x_2=\gamma(i+2) \rightarrow \cdots \rightarrow \gamma(n).$$ We thus have $|\Gamma_{x_0,x_1,x_2}| = A(x_0) B(x_2)$, where $A(x_0)$ is the number of oriented paths joining some $\gamma(0)$ such that $\mathcal{E}(\gamma(0))=\emptyset$ to $x_0$ and where $B(x_2)$ is defined similarly. $\square$

\section{General bounds on $H''(t)$} \label{sec:GeneralFormula}

In this section, we adapt the method used in~\cite{HillionContraction} to prove the convexity of the entropy along the contraction of a probability measure on $\mathbb{Z}$ to the more general framework of $W_{1,+}$-geodesics on a graph. We then apply this method in the cases where $G$ is the graph $\mathbb{Z}$ or a complete graph. We finally study the behaviour, along a $W_{1,+}$-geodesic, of the relative entropy with respect to a log-concave reference probability measure and discuss why the hypothesis of a uniform bound on the second derivative $H''(t)$ may not be by itself a sufficient condition for interseting functional inequalities to hold.

\subsection{Benamou-Brenier triples}

Let $G$ be a graph, endowed with the $W_1$-orientation with respect to a couple of probability distributions $f_0,f_1$ on $G$.

\begin{defi}
A Benamou-Brenier triple, or BB-triple, on $(G,\rightarrow)$, is a triple of positive functions $f,g,h$ defined respectively on $G$, $E(G)$ and $T(G)$ such that
\begin{equation}
\forall (x_0x_1x_2) \in T(G) \ , \ h(x_0x_1x_2)f(x_1) = g(x_0x_1)g(x_1x_2).
\end{equation}
\end{defi}

It is clear that, if a triple $(f_t,g_t,h_t)$ defines a $W_{1,+}$-geodesic, then for each $t \in [0,1]$, $(f_t,g_t,h_t)$ is a BB-triple. Other types of BB-triples will be considered in Section~\ref{sec:ProductSpaces}.

\begin{defi}
The functional $\mathcal{I}$ is defined for every BB-triple on $(G,\rightarrow)$ by
\begin{equation}\label{eq:defI}
\mathcal{I}(f,g,h) := \sum_{x \in G} \nabla_2 \cdot h(x) \log(f(x)) + \frac{(\nabla \cdot g(x))^2}{f(x)}.
\end{equation}
\end{defi}

\begin{prop}
Let us consider $(f_t,g_t,h_t)$ defining a $W_{1,+}$-geodesic on $G$. The entropy $H(t)$ of $f_t$ satisfies
\begin{equation}
H''(t) = \mathcal{I}(f_t,g_t,h_t).
\end{equation}
\end{prop}

\textbf{Proof: } This simply comes from the definition of the families $(g_t)_{t \in [0,1]}$ and $(h_t)_{t \in [0,1]}$:
\begin{equation}
\frac{\partial}{\partial t}f_t(x) = - \nabla \cdot  g_t(x) \ , \ \frac{\partial^2}{\partial t^2} f_t(x) = \nabla_2 \cdot h_t(x). \ \square
\end{equation}

\subsection{Integration by parts on $G$}

In order to obtain bounds on $H''(t)$, we first use integration by parts to transform the sum in~\eqref{eq:defI}:

\begin{prop} \label{prop:Itelescopic}
For any BB-triple $(f,g,h)$ we have
\begin{eqnarray*}
\mathcal{I}(f,g,h) &=& \sum_{x \in G} \left[ \sum_{x_1 \in \mathcal{F}(x)} \sum_{x_2 \in \mathcal{F}(x_1)} h(xx_1x_2) \log\left(\frac{f(x) h_t(x x_1 x_2)}{g(x x_1)^2}\right) \right] \\
&& + \sum_{x \in G} \left[ \sum_{x_{-1} \in \mathcal{E}(x)} \sum_{x_{-2} \in \mathcal{E}(x_{-1})} h(x_{-2}x_{-1}x) \log\left(\frac{f(x) h_t(x_{-2}x_{-1}x)}{g(x_{-1}x)^2}\right) \right] \\
&& + \sum_{x \in G} \frac{(\nabla \cdot g(x))^2}{f(x)}.
\end{eqnarray*}
\end{prop}

\textbf{Proof: } We add to the sum defining $\mathcal{I}(f,g,h)$ (see equation~\eqref{eq:defI}) the following telescopic sums
\begin{equation}
0 = \sum_{x \in G} \nabla_2 \cdot (h \log(h))(x) \ , \ 0= -2 \sum_{x \in G} \nabla \cdot \left(g \nabla \cdot  h\right)(x).
\end{equation}
The proposition is then proven by noticing that, $(f,g,h)$ being a BB-triple, we have
\begin{equation}
\forall x \in G \ , \  -2 \sum_{x_{-1} \in \mathcal{E}(x)} \sum_{x_1 \in \mathcal{F}(x)} h(x_{-1}xx_1) \log\left(\frac{h(x_{-1}xx_1) f(x)}{g(x_{-1}x)g(xx_1)} \right) =0. \  \square
\end{equation}

Combining Proposition~\ref{prop:Itelescopic} with the elementary inequality $\log(x) \geq 1-1/x$ allows us to obtain bounds on $\mathcal{I}(f,g,h)$:

\begin{prop}\label{prop:Iboundnaif}
For any triple $(f,g,h)$ we have
\begin{equation} \label{eq:Iboundnaif}
\mathcal{I}(f,g,h) \geq \sum_{(x_0x_1) \in E(G)} \frac{g(x_0x_1)^2}{f(x_0)}\left(1-\left|\mathcal{F}(x_1)\right| \right) + \frac{g(x_0x_1)^2}{f(x_1)}\left(1-\left|\mathcal{E}(x_0)\right| \right).
\end{equation}
\end{prop}

\textbf{Proof: } The inequality $\log(x) \geq 1-1/x$ implies
\begin{eqnarray*}
\mathcal{I}(f,g,h) &\geq& \sum_{x \in G} \left[ \sum_{x_1 \in \mathcal{F}(x)} \sum_{x_2 \in \mathcal{F}(x_1)} h(xx_1x_2) - \frac{g(x x_1)^2}{f(x)} \right] \\
&& + \sum_{x \in G} \left[ \sum_{x_{-1} \in \mathcal{E}(x)} \sum_{x_{-2} \in \mathcal{E}(x_{-1})} h(x_{-2}x_{-1}x) - \frac{g(x_{-1}x)^2}{f(x)} \right] \\
&& + \sum_{x \in G} \frac{(\nabla \cdot g(x))^2}{f(x)}.
\end{eqnarray*}

The following are obvious:
\begin{equation}
\sum_{x_2 \in \mathcal{F}(x_1)} \frac{g(x x_1)^2}{f(x)} = |\mathcal{F}(x_1)| \frac{g(x x_1)^2}{f(x)} \ , \  \sum_{x_{-2} \in \mathcal{E}(x_{-1})} \frac{g(x_{-1}x)^2}{f(x)} = |\mathcal{E}(x_{-1})| \frac{g(x_{-1}x)^2}{f(x)}.
\end{equation}

Moreover, we have:
\begin{eqnarray*}
\sum_{x \in G} \sum_{x_1 \in \mathcal{F}(x)} \sum_{x_2 \in \mathcal{F}(x_1)} h(xx_1x_2) &=& \sum_{(x_0x_1x_2) \in T(G)} h(x_0x_1x_2) \\
&=& \sum_{x \in G} \sum_{x_{-1} \in \mathcal{E}(x)} \sum_{x_1 \in \mathcal{F}(x)} \frac{g(x_{-1}x)g(xx_1)}{f(x)},
\end{eqnarray*}
and similarly:
\begin{equation}
\sum_{x \in G} \sum_{x_1 \in \mathcal{F}(x)} \sum_{x_2 \in \mathcal{F}(x_1)} h(xx_1x_2) = \sum_{x \in G} \sum_{x_{-1} \in \mathcal{E}(x)} \sum_{x_1 \in \mathcal{F}(x)} \frac{g(x_{-1}x)g(xx_1)}{f(x)}.
\end{equation}

Expanding $\sum_{x \in G} \frac{(\nabla \cdot g(x))^2}{f(x)}$ allows us to find similar terms:
\begin{eqnarray*}
\sum_{x \in G} \frac{(\nabla \cdot g(x))^2}{f(x)} &=& \sum_{x \in G} \frac{\left(\sum_{x_1 \in \mathcal{F}(x)} g(xx_1) - \sum_{x_{-1} \in \mathcal{E}(x)} g(x_{-1}x)\right)^2}{f(x)} \\
&\geq & \sum_{x \in G} \sum_{x_1 \in \mathcal{F}(x)} \frac{g(xx_1)^2}{f(x)}+\sum_{x \in G} \sum_{x_{-1} \in \mathcal{E}(x)} \frac{g(x_{-1}x)^2}{f(x)} \\
&& -2 \sum_{x \in G} \sum_{x_{-1} \in \mathcal{E}(x)} \sum_{x_1 \in \mathcal{F}(x)} \frac{g(x_{-1}x)g(xx_1)}{f(x)}.
\end{eqnarray*}
We used the fact that $g$ is non-negative to apply the inequality 
\begin{equation}
\left(\sum_{x_1 \in \mathcal{F}(x)} g(xx_1) \right)^2 \geq \sum_{x_1 \in \mathcal{F}(x)} g(xx_1)^2,
\end{equation}
which is far from being optimal, unless $|\mathcal{F}(x)|=0$ or $1$.

\medskip

Combining these estimations leads to 
\begin{eqnarray*}
\mathcal{I}(f,g,h) &\geq& \sum_{x \in G} \sum_{x_1 \in \mathcal{F}(x)} \left(1-|\mathcal{F}(x_1)|\right) \frac{g(x x_1)^2}{f(x)}\\&&+\sum_{x \in G} \sum_{x_{-1} \in \mathcal{E}(x)} \left(1-|\mathcal{E}(x_{-1})|\right)\frac{g(x_{-1}x)^2}{f(x)},
\end{eqnarray*}
which, up to a change of indices, is exactly inequality~\eqref{eq:Iboundnaif}. $\square$

\medskip

The bound obtained in Proposition~\ref{prop:Iboundnaif} is interesting in two fundamental cases:

\begin{cor} \label{cor:IboundCube}
Let $(f,g,h)$ be a BB-triple of functions on $(G,\rightarrow)$ where $G$ is the complete graph with $n$ points, $W_1$-oriented with respect to some couple $(f_0,f_1)$. We have
\begin{equation}
\mathcal{I}(f,g,h) \geq \sum_{(x_0x_1) \in E(G)} g(x_0x_1)^2 \left(\frac{1}{f(x_0)}+\frac{1}{f(x_1)} \right).
\end{equation}
\end{cor}

\textbf{Proof: } We apply Proposition~\ref{prop:Iboundnaif}, using the fact that, if $(x_0x_1) \in E(G)$, then the sets $\mathcal{E}(x_0)$ and $\mathcal{F}(x_1)$ are empty, or the equivalent fact that the set of oriented triple $T(G)$ is empty: indeed if there exists $(x_0x_1x_2) \in T(G)$ then, by Proposition~\ref{prop:OrientedGeodesic}, we have $d(x_0,x_2)=2$, which is a contradiction. $\square$

\begin{cor} \label{cor:IboundZ}
Let $(f,g,h)$ be a (finitely supported) BB-triple of functions on $(G,\rightarrow)$ where $G$ is the graph $\mathbb{Z}$, $W_1$-oriented with respect to some couple $(f_0,f_1)$. T $\mathcal{I}(f,g,h) \geq 0.$
\end{cor}

\textbf{Proof: } We use this time the fact that each vertex of $\mathbb{Z}$ has two neighbours, which implies that, for every $x \in \mathbb{Z}$, $|\mathcal{E}(x)|+|\mathcal{F}(x)| \leq 2$. In particular, if $(x_0x_1) \in E(G)$, then $\mathcal{E}(x_1)$ is non-empty (as it contains $x_0$), so $|\mathcal{F}(x_1)| \leq 1$. Similarly we have $|\mathcal{E}(x_0)| \leq 1$. Applying Proposition~\ref{prop:Iboundnaif} leads to the result. $\square$

\medskip

\textbf{Remark.} Corollary~\ref{cor:IboundZ} can be extended to the framework of cyclic graphs $\mathbb{Z}_r$ for $r \geq 2$, because in this case every vertex has also two neighbours.

\subsection{About the convexity of the relative entropy}

We have been so far interested in the behaviour of the Shannon entropy functional $H(f) := \sum_{x \in G} f(x) \log(f(x))$ along $W_{1,+}$ geodesics on $G$. However, the functional which is considered in Sturm-Lott-Villani theory are the relative entropy $H_\nu$ with respect to some reference probability measure $\nu$. In this paragraph, we present some results about the behaviour of $H_\nu$ along $W_{1,+}$-geodesics on graphs.

\begin{defi}
Let $\nu$ be a probability measure fully supported on $G$. The relative entropy $H_\nu(f)$ of a probability measure $f$ on $G$ is defined by 
\begin{equation*} 
H_\nu(f) := \sum_{x \in G} f(x) \log\left(\frac{f(x)}{\nu(x)}\right).
\end{equation*}
\end{defi}

\textbf{Remark.} Let $(f_t)_{t \in [0,1]}$ be a $W_{1,+}$-geodesic supported on a finite subset of vertices $A \subset G$. Let $\nu$ be the uniform probability distribution on $A$. Then the Shannon and relative entropies are linked by \begin{equation*} H_\nu(f_t) = H(f_t) + \log(|A|)\end{equation*} so the convexity of $t \mapsto H_\nu(f_t)$ is equivalent to the convexity of $t \mapsto H(f_t)$.

\medskip

As in the Riemannian case, it is interesting to consider log-concave reference measures:

\begin{prop}
We endow $G$ with a reference measure $\nu(x) := \exp(-V(x))$. We suppose that there exists $K>0$ such that, for every geodesic path of length $2$ $\gamma_0,\gamma_1,\gamma_2$ we have $$V(\gamma_0)-2 V(\gamma_1) + V(\gamma_2) \geq K.$$ Let $(f_t)$ be a $W_{1,+}$-geodesic, $H(t)$ be the Shannon entropy of $f_t$ and $H_\nu(t)$ its relative entropy. Then 
\begin{equation} \label{eq:EntroRelative}
H_\nu''(f_t) \geq  H''(t) + K W^2(f_0,f_1),
\end{equation} where $$W^2(f_0,f_1) := \sum_{(x_0x_1x_2) \in T(G)} h_t(x_0x_1x_2)$$ does not depend on $t$.
\end{prop}

\textbf{Proof:} We have 
\begin{equation*}
H_\nu(t) - H(t) = - \sum_{x \in G} f_t(x) \log(\nu(x)) = \sum_{x \in G} f_t(x) V(x), 
\end{equation*}
and by differentiating twice with respect to $t$ we have
\begin{equation*}
H_\nu''(t) = H''(t) + \sum_{x \in G} (\nabla_2 \cdot h_t)(x) V(x) = \sum_{(x_0x_1x_2) \in T(G)} h_t(x) (V(x_2)-2V(x_1)+V(x_0)),
\end{equation*}
which, by the convexity assumption made on $V$, proves equation~\eqref{eq:EntroRelative}. 

\medskip

In order to prove that $\sum_{(x_0x_1x_2) \in T(G)} h_t(x_0x_1x_2)$ does not depend on $t$, we first use the Benamou-Brenier condition~\eqref{eq:BBcondition} to write
\begin{equation*}
\frac{\partial}{\partial t} h_t(x_0x_1x_2) = - \sum_{x_{-1} \in \mathcal{E}(x_0)} \frac{g_t(x_{-1}x_0) g_t(x_0x_1) g_t(x_1x_2)}{f_t(x_0)f_t(x_1)}+ \sum_{x_3 \in \mathcal{F}(x_2)} \frac{g_t(x_0x_1)g_t(x_1x_2)g_t(x_2x_3)}{f_t(x_1)f_t(x_2)}.
\end{equation*}
A simple change of indices then show that
\begin{equation*}
\frac{\partial}{\partial t} \sum_{(x_0x_1x_2) \in T(G)} h_t(x_0x_1x_2) = \sum_{x_0 \rightarrow \cdots \rightarrow x_3 \in G} \frac{g_t(x_0x_1)g_t(x_1x_2)g_t(x_2x_3)}{f_t(x_1)f_t(x_2)}-\frac{g_t(x_0x_1)g_t(x_1x_2)g_t(x_2x_3)}{f_t(x_1)f_t(x_2)}=0,
\end{equation*}
so $\sum_{(x_0x_1x_2) \in T(G)} h_t(x_0x_1x_2)$ does not depend on $t$. $\square$

\medskip

\textbf{Remark.} One major difference with the continuous case is the fact that, although acting as the Wasserstein distance $W_2$, the quantity $W(f_0,f_1)$ does not define a distance on $\mathcal{P}(G)$. For instance, if $f_0$ and $f_1$ are two Dirac distributions at two adjacent vertices, we have $W(f_0,f_1)=0$. A different perspective on the same issue consists in writing $$W^2(f_0,f_1) = \sum_{x_1 \in G} f_t(x_1) V_{+,t}(x_1) V_{-,t}(x_1),$$ where $V_{+,t}(x_1) := \sum_{x_2 \in \mathcal{F}(x_1)} \frac{g_t(x_1x_2)}{f_t(x_1)}$ and $V_{-,t}(x_1) := \sum_{x_0 \in \mathcal{E}(x_1)} \frac{g_t(x_0x_1)}{f_t(x_1)}$ are the two velocity functions, which can be written $W^2 = \langle V_{+,t} , V_{-,t}\rangle$ for the scalar product with respect to $f_t$. This formula is the discrete analogue of the Benamou-Brenier formula~\eqref{eq:MKBenamou} for $p=2$, but in the continuous setting we have $W_2^2 = <v_t,v_t>= ||v_t||^2$ for the scalar product with respect to $f_t$. The fact that $V_{+,t} \neq V_{-,t}$ is a major obstacle to a generalization of the HWI inequality which holds for instance in the measured length space $(\mathbb{R}^d,\exp(-V(x))dx)$ (see~\cite{LottVillani} for a proof of this fact).

\section{Product of graphs} \label{sec:ProductSpaces}

Let $G_1$ and $G_2$ be two locally finite and connected graphs. In this section we study the behaviour of the entropy along $W_{1,+}$-geodesics defined on the product graph $G := G_1 \times G_2$ endowed with the usual product metric $$d_G((x_1,x_2),(y_1,y_2)) := d_{G_1}(x_1,y_1)+d_{G_2}(x_2,y_2).$$ 

\subsection{The $W_1$-orientation on a product graph}

The neighbours of a vertex $(x_1,x_2)$ in $G$ are the vertices $(x_1,y_2)$, where $d_{G_2}(x_2,y_2)=1$ and $(y_1,x_2)$ where $d_{G_1}(x_1,y_1)=1$. From this fact we easily deduce the following description of geodesic curves in $G$:

\begin{prop} \label{prop:GeodProduct}
Let $\gamma \in \Gamma(x,y)$ be a geodesic on $G$, where $(x,y)=((x_1,x_2),(y_1,y_2))$. There exist two geodesics $\gamma_1 \in \Gamma(x_1,y_1)$, $\gamma_2 \in \Gamma(x_2,y_2)$ defined respectively on $G_1$ and $G_2$, and an application $$\phi : \{0,\ldots, d(x,y)\} \rightarrow \{0,\ldots, d_1(x_1,y_1)\}$$ with $\phi(0)=0$, $\phi(d(x,y))=d_1(x_1,y_1)$ and $\phi(k+1)-\phi(k) \in \{0,1\}$, such that $$\gamma(k) = (\gamma_1(\phi(k)),\gamma_2(k-\phi(k))).$$
\end{prop}

In particular, the cardinality of $\Gamma(x,y)$ satisfies
\begin{equation}
|\Gamma(x,y)| = \binom{d(x,y)}{d(x_1,y_1)}|\Gamma(x_1,y_1)||\Gamma(x_2,y_2)|.
\end{equation}

If $f$ is a probability distribution on $G$, we denote by $f^{(1)}$, $f^{(2)}$ its marginals on $G_1$ and $G_2$. 
To a coupling $\pi$ between two distributions $f_0,f_1$, which can be seen as a probability measure on $$G \times G = (G_1 \times G_2) \times (G_1 \times G_2) = (G_1 \times G_1) \times (G_2 \times G_2),$$ we associate the marginal couplings $\pi^{(1)}$ on $G_1 \times G_1$ between $f_0^{(1)}$ and $f_1^{(1)}$ and $\pi_2$ on $G_2 \times G_2$ between $f_0^{(2)}$ and $f_1^{(2)}$.

\medskip

We then describe the $W_{1,+}$-orientation on $G$ with respect to a couple of measures $f_0,f_1$.

\begin{prop} \label{prop:W1OrientProduct}
Let $f_0,f_1 \in \mathcal{P}(G)$. For $i=1,2$ we define
\begin{equation}
\mathcal{E}_i(x^{(i)}) := \{y^{(i)} \in G_i \ : \ y^{(i)} \rightarrow x^{(i)} \} \ , \ \mathcal{F}_i(x^{(i)}) := \{z^{(i)} \in G_i \ : \ x^{(i)} \rightarrow z^{(i)} \}
\end{equation}
for the $W_1$ orientation on $G_i$ between $f_0^{(i)}$ and $f_1^{(i)}$. The $W_1$-orientation between $f_0$ and $f_1$ is then described by
\begin{eqnarray}
\mathcal{E}(x) &=&  \left( \bigcup_{y^{(2)} \in \mathcal{E}_2(x^{(2)})} (x^{(1)},y^{(2)}) \right) \bigcup \left( \bigcup_{y^{(1)} \in \mathcal{E}_2(x^{(1)})} (y^{(1)},x^{(2)}) \right)\\
&=:& \mathcal{E}_1(x) \cup \mathcal{E}_2(x), \\
\mathcal{F}(x) &=&  \left( \bigcup_{y^{(2)} \in \mathcal{F}_2(x^{(2)})} (x^{(1)},y^{(2)}) \right) \bigcup \left( \bigcup_{y^{(1)} \in \mathcal{F}_2(x^{(1)})} (y^{(1)},x^{(2)}) \right)\\
&=:& \mathcal{F}_1(x) \cup \mathcal{F}_2(x).
\end{eqnarray}  
\end{prop}

\textbf{Proof: } Let $\pi \in \Pi(f_0,f_1)$ be a coupling between $f_0$ and $f_1$. We have
\begin{eqnarray*}
 I_1(\pi) &=& \sum_{(x_1,x_2),(y_1,y_2) \in G \times G} d((x_1,x_2),(y_1,y_2)) \pi((x_1,x_2),(y_1,y_2)) \\
 &=& \sum_{(x_1,y_1) \in G_1 \times G_1} \sum_{(x_2,y_2) \in G_2\times G_2} d_1(x_1,y_1)+d_2(x_2,y_2) \pi((x_1,x_2),(y_1,y_2)) \\
 &=& \sum_{(x_1,y_1) \in G_1 \times G_1} d_1(x_1,y_1) \pi^{(1)}(x_1,y_1) + \sum_{(x_2,y_2) \in G_2\times G_2}d_2(x_2,y_2) \pi^ {(2)}(x_2,y_2) \\
 &=& I_1(\pi^{(1)})+ I_1(\pi^{(2)}),
\end{eqnarray*}
which proves that $\pi$ is $W_1$-optimal between $f_0$ and $f_1$ (for the distance $d_G$) if and only if its marginals $\pi^{(1)}$, $\pi^{(2)}$ are $W_1$-optimal between $f_0^{(1)}$ and $f_1^{(1)}$, resp $f_0^{(2)}$ and $f_1^{(2)}$ for the distance $d_{G_1}$, resp. $d_{G_2}$.

\medskip

We now fix a $W_1$-optimal coupling $\pi \in \Pi_1(f_0,f_1)$. Let $x=(x_1,x_2)$ and $y=(y_1,y_2)$ be two vertices of $G$ such that $\pi(x,y)>0$. We then have $\pi^{(1)}(x_1,y_1)>0$ and $\pi^{(2)}(x_2,y_2)>0$ for the marginal couplings, which are also $W_1$-optimal. 

\medskip

Let $\gamma \in \Gamma_G(x,y)$, and $\gamma_1 \in \Gamma_{G_1}(x_1,y_1)$, $\gamma_2 \in \Gamma_{G_2}(x_2,y_2)$, $\phi : \{0,\ldots, d(x,y)\} \rightarrow \{0,\ldots, d_1(x_1,y_1)\}$ be associated to $\gamma$ by Proposition~\ref{prop:GeodProduct}. For $k \in \{0,\ldots, d(x,y)-1\}$, we have $$\gamma(k) = (\gamma_1(\phi(k)),\gamma_2(k-\phi(k))) \ , \ \gamma(k+1) = (\gamma_1(\phi(k)),\gamma_2(k-\phi(k+1)+1)).$$ If $\phi(k+1)=\phi(k)+1$, resp. $\phi(k+1)=\phi(k)$, then $\gamma(k+1) \in \mathcal{E}_1(\gamma(k))$, resp. $\gamma(k+1) \in \mathcal{E}_2(\gamma(k))$.

\medskip

Conversely, let us consider a vertex $x=(x_1,x_2) \in G$. We suppose that $x=\gamma(k)$ for some geodesic $\gamma$ of length $n$ such that $\pi(\gamma(0),\gamma(n))>0$ for a $W_1$-optimal coupling $\pi \in \Pi_1(f_0,f_1)$. We denote by $\gamma^{(1)}, \gamma^{(2)}$ the projections of $\gamma$, as defined in Proposition~\ref{prop:GeodProduct}, $n_1$ and $n_2$ their respective lengths, and $\pi^{(1)}$, $\pi^{(2)}$ the marginals of $\pi$. Let $y \in \mathcal{E}_1(x)$. We have $y=(y_1,x_2)$ with $y_1 \in \mathcal{E}_1(x_1)$. There exists a $W_1$-optimal coupling $\tilde{\pi}^{(1)} \in \Pi_1(f_0,f_1)$ and a geodesic $\tilde{\gamma}_1$ on $G_1$, of length $\tilde{n}_1$, such that $\tilde{\gamma}_1(k_1)=x_1$, $\tilde{\gamma}_1(k_1+1)=y_1$ and $\tilde{\pi}^{(1)}(\tilde{\gamma}_1(0),\tilde{\gamma}_1(\tilde{n}_1)) > 0$. Let $\tilde{\pi}$ be any coupling between $f_0$ and $f_1$ having $\tilde{\pi}^{(1)}$ and $\pi^{(2)}$ as marginals and $\gamma$ be a geodesic of $G$ having $\tilde{\gamma_1}$ and $\gamma_2$ as 
projections. Then there exists some $k$ for which $\tilde{\gamma}(k)=x$, $\tilde{\gamma}(k+1)=y$. Furthermore $\tilde{\pi}$ is $W_1$-optimal between $f_0$ and $f_1$ and $\tilde{\pi}(\gamma(0),\gamma(\tilde{n}_1+n_2)) > 0$, which proves that $y \in \mathcal{E}(x)$.

\medskip

We can prove similarly that, if $y \in \mathcal{E}_2(x)$ then $y \in \mathcal{E}(x)$, which finishes the proof. $\square$

\medskip

An immediate consequence of Proposition~\ref{prop:W1OrientProduct} is a decomposition of the divergence operator:

\begin{prop} \label{prop:divdecomp}
The divergence $\nabla \cdot  g$ of a function $g : E(G) \rightarrow \mathbb{R}$ can be written $\nabla \cdot  g= \nabla^{(1)} \cdot g+\nabla^{(2)} \cdot g$ where
\begin{equation}
\nabla^{(i)} \cdot g(x_1) := \sum_{x_2 \in \mathcal{F}_i(x_1)} g(x_1 x_2) - \sum_{x_0 \in \mathcal{E}_i(x_1)} g(x_0x_1).
\end{equation}
Similarly, the second order divergence operator of a function $h : T(G) \rightarrow \mathbb{R}$ can be written
\begin{equation}
\nabla_2 \cdot h = \nabla_2^{(11)} \cdot h +\nabla_2^{(12)} \cdot h+\nabla_2^{(21)} \cdot h+\nabla_2^{(22)} \cdot h,
\end{equation}
with $\nabla_2^{(ij)} := \nabla^{(i)} \circ \nabla^{(j)}$.
\end{prop}

The structure of the oriented graph $(G_1 \times G_2)$ is better understood by introducing oriented product squares:

\begin{defi}
An oriented product square of $G$ is a 4-uple of vertices $(x_0,x_1,x_1',x_2) \in G^4$ such that $x_1 \in \mathcal{F}_1(x_0)$, $x_1' \in \mathcal{F}_2(x_0)$, $x_2 \in \mathcal{F}_2(x_1)$ and $x_2 \in \mathcal{F}_1(x_1')$. We denote by $S(G)$ the set of oriented product squares of $G$.
\end{defi}

\begin{prop}\label{prop:OrientedSquares}
Let $x_0 \in G$. The following sets are all in bijection: 
\begin{itemize}
\item $\mathcal{A}_1$:= $\mathcal{F}_1(x_0) \times \mathcal{F}_2(x_0)$.
\item $\mathcal{A}_2$:= $\{ x_2 \in G \  : \ \exists x_1,x_1' \in G\times G \ , \ (x_0,x_1,x_1',x_2) \in S(G) \}$.
\item $\mathcal{A}_3$:= $\{ (x_0x_1x_2) \in T(G) \ : x_1 \in \mathcal{F}_1(x_0) , x_2 \in \mathcal{F}_2(x_1) \ \}$.
\item $\mathcal{A}_4$:= $\{ (x_0x_1'x_2) \in T(G) \ : x_1' \in \mathcal{F}_2(x_0) , x_2 \in \mathcal{F}_1(x_1') \ \}$.
\end{itemize}
\end{prop}

\textbf{Proof:} Let us fix $x_1 \in \mathcal{F}_1(x_0)$ and $x_1' \in \mathcal{F}_2(x_0)$. We write $x_0=(x_0^{(1)},x_0^{(2)})$ in $G_1 \times G_2$. There exist a unique $x_1^{(1)} \in \mathcal{F}_{G_1}(x_0^{(1)})$ and a unique $x_1'^{(2)} \in \mathcal{F}_{G_{2}}(x_0^{(2)})$ such that $x_1 = (x_1^{(1)},x_0^{(2)})$ and $x_1'=(x_0^{(1)},x_1'^{(2)})$ in $G_1 \times G_2$. We then set $x_2 := (x_1^{(1)},x_1'^{(2)})$ and it is easy to see that $(x_0,x_1,x_1',x_2) \in S(G)$. $\square$

\medskip

Proposition~\ref{prop:OrientedSquares} shows that an oriented square $(x_0x_1x_1'x_2)$ is uniquely determined by the couple $x_0,x_2$. We will use the notation $(x_0x_2) \in S(G)$ to denote such squares. We will also denote the two midpoints $x_1,x_1'$ respectively by $m_1(x_0x_2)$ and $m_2(x_0x_2)$.

\medskip

Let $(f_t)$ be a $W_{1,+}$-geodesic on $G$. There exist two families of functions $(g_t)$ and $(h_t)$, defined respectively on $E(G)$ and $T(G)$, such that $\frac{\partial}{\partial t}f = - \nabla \cdot  g$, $\frac{\partial}{\partial t}g = - \nabla \cdot  h$ and satisfying $$ \forall (x_0x_1x_2) \in T(G) \ , \ f_t(x_1)h_t(x_0x_1x_2) = g_t(x_0x_1)g_t(x_1x_2).$$
Given a vertex $x^{(2)}$, we now define, for $(x_0^{(1)}x_1^{(1)}x_2^{(1)}) \in T(G_1)$, the functions
\begin{eqnarray} \label{eq:BBprojection}
f_{t,x^{(2)}}(x_0^{(1)}) &:=& f_t(x_0^{(1)},x^{(2)}), \\
g_{t,x^{(2)}}(x_0^{(1)}x_1^{(1)}) &:=& g_t((x_0^{(1)},x^{(2)})(x_1^{(1)},x^{(2)})), \\
h_{t,x^{(2)}}(x_0^{(1)}x_1^{(1)}x_2^{(1)}) &:= & h_t((x_0^{(1)},x^{(2)})(x_1^{(1)},x^{(2)})(x_2^{(1)},x^{(2)})).
\end{eqnarray}
The triple of functions $(f_{t,x^{(2)}},g_{t,x^{(2)}},h_{t,x^{(2)}})$ is then a BB-triple on $G_1$. Given $x^{(1)} \in G_1$, we define similarly the BB-triples of functions $(f_{t,x^{(1)}},g_{t,x_{(1)}},h_{t,x_{(1)}})$ on $G_2$.

\medskip

The divergence of $g_{t,x^{(2)}} : E(G_1) \rightarrow \mathbb{R}$ satisfies the relation
\begin{equation}
(\nabla \cdot  g_{t,x^{(2)}})(x^{(1)}) = (\nabla^{(1)} \cdot g_t)(x^{(1)},x^{(2)}).
\end{equation}
The second order divergence $h_{t,x^{(2)}} : T(G_1) \rightarrow \mathbb{R}$ satisfies
\begin{equation}
(\nabla_2 \cdot h_{t,x^{(2)}})(x^{(1)}) = (\nabla^{(11)} \cdot h_t)(x^{(1)},x^{(2)}).
\end{equation}

\subsection{A tensorization result}

We are now able to state the tensorization theorem:

\begin{thm} \label{th:EntroTensorization}
Let $(f_t,g_t,h_t)$ be a canonical $W_{1,+}$-geodesic on $G$ and $H(t)$ denote the entropy of $f_t$. Then:
\begin{equation}
H''(t)  \geq \sum_{x^{(2)} \in G_2} \mathcal{I}(f_{t,x^{(2)}},g_{t,x^{(2)}},h_{t,x^{(2)}}) + \sum_{x^{(1)} \in G_1} \mathcal{I}(f_{t,x^{(1)}},g_{t,x^{(1)}},h_{t,x^{(1)}}).
\end{equation}
\end{thm}

\textbf{Proof: }We apply Proposition~\ref{prop:divdecomp}:
\begin{eqnarray*}
\sum_{x^{(2)} \in G_2} \mathcal{I}(f_{t,x^{(2)}},g_{t,x^{(2)}},h_{t,x^{(2)}}) &=& \sum_{x^{(2)} \in G_2} \left( \sum_{x^{(1)} \in G_1} \nabla_2 \cdot h_{t,x^{(2)}}(x^{(1)}) \log(f_{t,x^{(2)}}(x^{(1)})) \right) \\
&& + \sum_{x^{(2)} \in G_2} \left( \sum_{x^{(1)} \in G_1} \frac{(\nabla \cdot  g_{t,x^{(2)}}(x^{(1)}))^2}{f_{t,x^{(2)}}(x^{(1)})} \right) \\ 
&=& \sum_{x \in G} \nabla_2^{(11)} \cdot h_t(x) \log(f_t(x)) + \frac{(\nabla^{(1)} \cdot g_t(x))^2}{f_t(x)}.
\end{eqnarray*}
Similarly, 
\begin{equation}
\sum_{x^{(1)} \in G_1} \mathcal{I}(f_{t,x^{(1)}},g_{t,x^{(1)}},h_{t,x^{(1)}}) = \sum_{x \in G} \nabla_2^{(22)} \cdot h_t(x) \log(f_t(x)) + \frac{(\nabla^{(2)} \cdot g_t(x))^2}{f_t(x)}.
\end{equation}
To prove Theorem~\ref{th:EntroTensorization}, it thus suffices to show the inequality
\begin{equation}\label{eq:CrossTerms}
\sum_{x \in G} (\nabla_2^{(12)}+\nabla_2^{(21)}) \cdot h_t(x) \log(f_t(x)) + 2 \frac{\nabla^{(1)} g_t(x) \nabla^{(2)} \cdot g_t(x)}{f_t(x)} \geq 0.
\end{equation}
By considering the telescopic sums $$\sum_{x \in G} \nabla^{(12)} \cdot  h_t \log(h_t)(x)=0 \ , \ -2\sum_{x\in G} \nabla^{(1)} \cdot (\nabla^{(2)} \cdot h_t \log(g_t))=0, $$ we prove, as in Proposition~\ref{prop:Itelescopic}, that
\begin{eqnarray*}
\sum_{x \in G} \nabla_2^{(12)} \cdot h_t(x) \log(f_t(x)) &=& \sum_{(x_0x_1x_2) \in T^{(12)}(G)} h(x_0x_1x_2) \log\left(\frac{f_t(x_0)h_t(x_0x_1x_2)}{g_t(x_0x_1)^2}\right)\\&&+\sum_{(x_0x_1x_2) \in T^{(12)}(G)} h(x_0x_1x_2)\log\left(\frac{f_t(x_2)h_t(x_0x_1x_2)}{g_t(x_1x_2)^2}\right),
\end{eqnarray*}
where $T^{(12)}(G)$ is the set of oriented triples $(x_0x_1x_2) \in T(G)$ such that $x_0 \in \mathcal{E}_1(x_1)$ and $x_1 \in \mathcal{E}_2(x_2)$. We now use the bijection between $T^{(12)}(G)$ and $S(G)$, proven in Proposition~\ref{prop:OrientedSquares}, and the fact that $h(x_0x_1x_2)$ does not depend on $x_1$, which comes from the assumption that $(f_t)$ is canonical and from  Proposition~\ref{prop:CanonicalTriple}, to write:
\begin{equation*}
\sum_{x \in G} \nabla_2^{(12)} \cdot h_t(x) \log(f_t(x)) = \sum_{(x_0x_2) \in S(G)} h(x_0x_2) \left( \log \left(\frac{f(x_0)h(x_0x_2)}{g(x_0m_1(x_0,x_2))^2 } \right)+\log \left(\frac{f(x_2)h(x_0x_2)}{g(m_1(x_0,x_2)x_2)^2 } \right)\right).
\end{equation*}
Similarly, we have:
\begin{equation*}
\sum_{x \in G} \nabla_2^{(21)} \cdot h_t(x) \log(f_t(x)) = \sum_{(x_0x_2) \in S(G)} h(x_0x_2) \left( \log \left(\frac{f(x_0)h(x_0x_2)}{g(x_0m_2(x_0,x_2))^2 } \right)+\log \left(\frac{f(x_2)h(x_0x_2)}{g(m_2(x_0,x_2)x_2)^2 } \right)\right).
\end{equation*}
Adding both equations and using the inequality $\log(x) \geq 1-1/x$ gives:
\begin{eqnarray*}
\sum_{x \in G} (\nabla_2^{(12)}+\nabla_2^{(21)}) \cdot h_t(x) \log(f_t(x)) &=& 2 \sum_{(x_0x_2) \in S(G)} h(x_0x_2) \log \left(\frac{f(x_0)h(x_0x_2)}{g(x_0m_1(x_0,x_2))g(x_0m_2(x_0,x_2)) } \right) \\&&+ 2 \sum_{(x_0x_2) \in S(G)} h(x_0x_2) \log \left(\frac{f(x_2)h(x_0x_2)}{g(m_1(x_0,x_2)x_2)g(m_2(x_0,x_2)x_2) } \right) \\
&\geq& 4 \sum_{(x_0x_2) \in S(G)} h(x_0x_2) \\ && - 2 \sum_{(x_0x_2) \in S(G)} \frac{g(x_0m_1(x_0,x_2))g(x_0m_2(x_0,x_2))}{f_t(x_0)} \\ && - 2 \sum_{(x_0x_2) \in S(G)} \frac{g(m_1(x_0,x_2)x_2)g(m_2(x_0,x_2)x_2)}{f_t(x_2)}.
\end{eqnarray*}
We use again the bijection in Proposition~\ref{prop:OrientedSquares} to write 
\begin{eqnarray*}
\sum_{(x_0x_2) \in S(G)} \frac{g(x_0m_1(x_0,x_2))g(x_0m_2(x_0,x_2)}{f_t(x_0)} &=& \sum_{x_0 \in G} \left( \sum_{(x_1,x_1') \in \mathcal{F}_1(x_0) \times \mathcal{F}_2(x_0)} \frac{g(x_0x_1)g(x_0x_1')}{f(x_0)} \right)\\
&=& \sum_{x_0 \in G}  \frac{\sum_{x_1 \in \mathcal{F}_1(x_0)} g(x_0x_1) \cdot \sum_{x_1' \in \mathcal{F}_2(x_0)} g(x_0x_1')}{f(x_0)} 
\end{eqnarray*}
and
\begin{eqnarray*}
\sum_{(x_0x_2) \in S(G)} \frac{g(m_1(x_0,x_2)x_2)g(m_2(x_0,x_2)x_2)}{f_t(x_2)} = \sum_{x_0 \in G}  \frac{\sum_{x_{-1} \in \mathcal{E}_1(x_0)} g(x_{-1}x_0) \cdot \sum_{x_{-1}' \in \mathcal{E}_2(x_0)} g(x_{-1}'x_0)}{f(x_0)}. 
\end{eqnarray*}
We also have:
\begin{eqnarray*}
\sum_{(x_0x_2) \in S(G)} h(x_0x_2) &=& \sum_{(x_{-1}x_0x_1) \in T^{(12)}(G)} h(x_{-1}x_0x_1) \\
&=& \sum_{x_0 \in G} \sum_{x_{-1} \in \mathcal{E}_1(x_0) \ , \ x_1 \in \mathcal{F}_2(x_0)} \frac{g(x_{-1}x_0) g(x_0x_1)}{f(x_0)} \\
&=& \sum_{x_0 \in G}  \frac{\sum_{x_{-1} \in \mathcal{E}_1(x_0)} g(x_{-1}x_0) \cdot \sum_{x_{1} \in \mathcal{F}_2(x_0)} g(x_{1}x_0)}{f(x_0)},
\end{eqnarray*}
and:
\begin{eqnarray*}
 \sum_{(x_0x_2) \in S(G)} h(x_0x_2) &=& \sum_{(x_{-1}x_0x_1) \in T^{(21)}(G)} h(x_{-1}x_0x_1) \\
 &=& \sum_{x_0 \in G}  \frac{\sum_{x_{-1} \in \mathcal{E}_2(x_0)} g(x_{-1}x_0) \cdot \sum_{x_{1} \in \mathcal{F}_1(x_0)} g(x_{1}x_0)}{f(x_0)}.
\end{eqnarray*}
Adding the last four identities gives:
\begin{eqnarray*}
\sum_{x \in G} (\nabla_2^{(12)}+\nabla_2^{(21)}) \cdot h_t(x) \log(f_t(x)) & \geq & -\sum_{x_0 \in G} \frac{\nabla^{(1)} \cdot g(x_0) \nabla^{(2)} \cdot g(x_0)}{f(x_0)},
\end{eqnarray*}
which is exactly the inequality~\eqref{eq:CrossTerms} we wanted to obtain. $\square$

\subsection{Examples}

The tensorization Theorem~\ref{th:EntroTensorization} is generalized to products of more than two graphs: let $G = G_1 \times \cdots \times G_p$. For $i = 1,\ldots,p$, we denote by $\hat{G}_i$ the product $G_1 \times \cdots G_p$, where $G_i$ is omitted. Given some vertex $\hat{x} \in \hat{G}_i$ and a BB-triple $(f_t,g_t,h_t)$ on $G$, we define a BB-triple $(f_{t,\hat{x}} ,g_{t,\hat{x}},h_{t,\hat{x}})$ as in equation~\eqref{eq:BBprojection}. We then have:

\begin{cor} \label{cor:EntroTensoMultiple}
Let $(f_t,g_t,h_t)$ be a BB-triple on $G$. The entropy $H(t)$ of $f_t$ satisfies 
\begin{equation}
H''(t) \geq \sum_{i=1}^p \sum_{\hat{x} \in \hat{G}_i} \mathcal{I}(f_{t,\hat{x}} ,g_{t,\hat{x}},h_{t,\hat{x}}).
\end{equation}
\end{cor}

Applying Corollary~\ref{cor:EntroTensoMultiple} to the examples studied in Section 3 allows us to obtain interesting bounds on the second derivative $H''(t)$ in other important cases:

\begin{prop} \label{prop:EntroTensoExemples}
Theorem~\ref{th:EntroTensorization} can be applied in the following fundamental examples:
\begin{itemize}
\item The entropy $H(t)$ along a $W_{1,+}$-geodesic  $(f_t)_{t \in [0,1]}$ on $\mathbb{Z}^n$ is a convex function of $t$.
\item Let $(f_t,g_t,h_t)$ be a $W_{1,+}$-geodesic on the cube $\mathbb{Z}_2^n$. Then 
\begin{equation}
\mathcal{I}(f_t,g_t,h_t) = \sum_{(x_0x_1) \in E(G)} g_t(x_0x_1)^2 \left(\frac{1}{f_t(x_0)} + \frac{1}{f_t(x_1)} \right) \geq 0.
\end{equation}
\end{itemize}
\end{prop}

\textbf{Proof:} The first point follows directly from Corollary~\ref{cor:IboundZ}. To prove the second point, we notice that the cube is described by the product $G_1 \times \cdots \times G_n$ where each $G_i$ is the two-point graph $\mathbb{Z}_2$. Each $\hat{G}_i$ is isometric to the $n-1$-dimensional cube. To each $\hat{x} \in \hat{G}_i$, we associate two vertices $\hat{x}_0,\hat{x}_1 \in G$ by setting the $i$-th coordiante to $0$ or $1$. If $\hat{x}_0 \rightarrow \hat{x}_1$ in $G$, we define $g_t(\hat{x}) := g_t(\hat{x}_0 \hat{x}_1)$. If $\hat{x}_1 \rightarrow \hat{x}_0$ we define $g_t(\hat{x}) := g_t(\hat{x}_1 \hat{x}_0)$. Finally if the edge $(\hat{x}_0 \hat{x}_1)$ is not oriented in $G$ we set $g_t(\hat{x}):=0$. In any case, we have, by Corollary~\ref{cor:IboundCube}, $$\mathcal{I}(f_{t,\hat{x}} ,g_{t,\hat{x}},h_{t,\hat{x}}) = g_t(\hat{x})^2 \left(\frac{1}{f_t(\hat{x}_0)}+\frac{1}{f_t(\hat{x}_1)} \right).$$
A (non-ordered) edge $(x_0 x_1)$ of $G$ is described in the following way: $x_0$ and $x_1$ differ by exactly one coordinate. In other terms, there is a bijection between the set of edges of $G$ and the disjoint union $\bigcup_{i=1}^p \hat{G}_i$. We can then write
\begin{equation}
\sum_{i=1}^n \sum_{\hat{x} \in \hat{G}_i}  \mathcal{I}(f_{t,\hat{x}} ,g_{t,\hat{x}},h_{t,\hat{x}}) = \sum_{(x_0x_1) \in (E(G),\rightarrow)} g_t(x_0x_1)^2 \left(\frac{1}{f_t(x_0)}+\frac{1}{f(x_1)} \right),
\end{equation}
which is what we wanted. $\square$

\medskip

These two examples can be seen as particular cases of a more general theorem:

\begin{thm}\label{th:GroupCurvature}
Let $G$ be the Cayley graph of a finitely generated abelian group, with a set of generators $T=(\tau_1,\ldots, \tau_q)$. Let $(f_t)$ be a $W_{1,+}$-interpolation on $G$ and $H(t)$ the entropy of $f_t$. Then :
\begin{equation} \label{eq:IboundGroup}
H''(t) \geq  \sum_{(x_0 x_1) \in \tilde{E}(G)} g_t(x_0x_1)^2 \left(\frac{1}{f_t(x_0)}+\frac{1}{f(x_1)} \right),
\end{equation}
where $\tilde{E}(G)$ is the subset of oriented edges $(x_0 \rightarrow x_1) \in E(G)$ such that $x_1 = \tau_i x_0$ for some generator $\tau_i \in T$ such that $\tau_i^2 = id$.
\end{thm}

\textbf{Proof: } Theorem~\ref{th:GroupCurvature} can be proven with the help of Theorem~\ref{th:EntroTensorization}. Indeed, any finitely generated abelian group is isomorphic to the direct product
\begin{equation*}
\mathbb{Z}^n \times \mathbb{Z}_2^{n_2} \times \cdots \times \mathbb{Z}_p^{n_p} \times \cdots, 
\end{equation*}
where all but a finite number of coefficients $n_p$ are equal to $0$. As we have proven that $\mathcal{I}(f,g,h) \geq 0$ for any BB-triple on $\mathbb{Z}_p$ or on $\mathbb{Z}$, a direct application of Theorem~\ref{th:EntroTensorization} gives that $H''(t) \geq 0$. The more precise bound given in equation~\eqref{eq:IboundGroup} is proven as in the second point of Proposition~\ref{prop:EntroTensoExemples}. $\square$

\appendix

\section{Appendix: further results on $W_{1,+}$-geodesics on $\mathbb{Z}$.}

\subsection{Renyi entropy along $W_{1,+}$-geodesics on $\mathbb{Z}$.}

In this appendix we prove that, along a $W_{1,+}$-geodesic on $\mathbb{Z}$, not only the relative entropy is convex, but also a larger class of functionals belonging to the family of Renyi entropies: given a probability distribution $(f(k))_{k \in \mathbb{Z}}$ and a parameter $0<p<1$, we set \begin{equation*} H_p(f) := -\sum_{k \in \mathbb{Z}} f(k)^p .\end{equation*} The relative entropy $H(f) := \sum_k f(k) \log(f(k))$ can be seen as a limit case of Renyi entropy as the parameter $p \rightarrow 1$, in the sense that \begin{equation*} H_p(f) = -1 + (1-p) H(f) + o((1-p)^2). \end{equation*}

We then have:
\begin{thm}\label{th:RenyiZ}
Let $(f_t)_{t \in [0,1]}$ be a $W_{1,+}$ geodesic on $\mathbb{Z}$. Then $t \mapsto H_p(f_t)$ is convex.
\end{thm}

In order to have simpler notations, we are going to prove Theorem~\ref{th:RenyiZ} under the additional assumption that $f_0$ is stochastically dominated by $f_1$ (see also Theorem~\ref{th:EntroBinoW2}). Under this assumption, the $W_1$-orientation on $\mathbb{Z}$ is simply described by orienting the edge $(k,k+1)$ by $k \rightarrow k+1$. If $g : E(G) \rightarrow \mathbb{R}$ is a function defined on oriented edges, we can then simply write $g(k)$ instead of $g(k,k+1)$ and the divergence operator $(\nabla \cdot g)(k) = (g(k)-g(k-1)$ can be seen as the left derivative of $g$. We will denote $\nabla g(k) := (\nabla \cdot g)(k)$. Similarly, if $k \rightarrow k+1 \rightarrow k+2$ is an oriented triple, we will write $h(k)$ instead of $h(k,k+1,k+2)$ and $(\nabla_2 \cdot h)(k)$ will be the twice left derivative $\nabla_2 h(k) := h(k)-2h(k-1)+h(k-2)$. With these notations, the Benamou-Brenier condition~\eqref{eq:BBcondition} is written \begin{equation} \label{eq:BBCondZ} h_t(k-1) f_t(k) = g_t(k-1) vg_t(k) . \end{equation}

\medskip

The proof of Theorem~\ref{th:RenyiZ} is based on two technical lemmas:

\begin{lem}
For every triple of non negative numbers $f,g,h$ we have
\begin{equation}\label{eq:HolderLike}
h f^{p-1} \leq \frac{1}{2-p} h^{2-p}g^{2p-2}+ \frac{1-p}{2-p} g^2 f^{p-2}.
\end{equation}
\end{lem}

\textbf{Proof: } The convexity of the exponential function implies the inequality:
\begin{equation}
\forall a,b >0 \ , \ \forall \alpha,\beta >0 \ s.t. \ \frac{1}{\alpha}+\frac{1}{\beta}=1 \ , \ ab \leq \frac{a^\alpha}{\alpha}+ \frac{b^\beta}{\beta},
\end{equation}
Setting $$\alpha := 2-p \ , \ \beta := \frac{2-p}{1-p} \ , \ a:= h g^{-2\frac{1-p}{2-p}} \ , \ b:= g^{2\frac{1-p}{2-p}} f^{p-1},$$ we obtain \eqref{eq:HolderLike} as wanted. $\square$

\begin{lem}
For every $x \geq 0$ we have
\begin{equation}
\psi(x) := (1-p)(x-1)^2 - \frac{x^p-1}{2-p}-\frac{1}{2-p}x^{2-p}-\frac{1-p}{2-p}x^2+2x-1 \geq 0.
\end{equation}
\end{lem}

\textbf{Proof: } Let us compute the first derivatives of $h$:
$$\psi'(x) = 2(1-p)(x-1) -\frac{p}{2-p} x^{p-1} -x^{1-p}-2\frac{1-p}{2-p}x+2,$$
$$\psi''(x) = 2(1-p) + \frac{p(1-p)}{(2-p)} x^{p-2}-(1-p)x^{-p}-2 \frac{1-p}{2-p},$$
$$\psi^{(3)}(x) = p(1-p) (x^{-p-1}-x^{p-3}).$$
As $p<1$, we have $-p-1>p-3$, so $\psi^{(3)}(x)$ is negative for $0 \leq x \leq 1$ and positive for $x\geq 1$. This means that $\psi''(x) \geq \psi''(1)=0$, so $\psi$ is convex. As we have $\psi(1)=0$ and $\psi'(1)=0$, we deduce that $\psi(x) \geq 0$ for every $x \geq 0$. $\square$

\medskip

\textbf{Proof of Theorem~\ref{th:RenyiZ}: } The second derivative of $H_p(f_t)$ satisfies:
\begin{equation*}
-H_p''(t) = p \sum_{k \in \mathbb{Z}} \nabla_2 h_t(k) f_t(k)^{p-1} + p(p-1) \sum_{k \in \mathbb{Z}} (\nabla g_t(k))^2 f_t(k)^{p-2}.
\end{equation*}
As $p \in (0,1)$, we have $p(1-p) < 0$. To prove the convexity of $t \mapsto H_p(t)$, we first apply twice the inequality \eqref{eq:HolderLike} with $(f,g,h)=(f_t(k),g_t(k),h_t(k))$ and $(f,g,h)=(f_t(k),g_t(k-1),h_t(k-2))$ and then change indices to write
\begin{eqnarray*}
\sum_{k \geq 0} \nabla_2(h_t(k)) f_t(k)^{p-1} &\leq& \sum_{k \geq 0} \frac{1}{2-p} h_t(k)^{2-p} g_t(k)^{2p-2}+ \frac{1-p}{2-p} g_t(k)^2 f_t(k)^{p-2}\\
&& -2h_t(k-1) f_t(k)^{p-1} \\ 
&&+ \frac{1}{2-p} h_t(k-2)^{2-p} g_t(k-1)^{2p-2} + \frac{1-p}{2-p} g_t(k-1)^2 f_t(k)^{p-2} \\
&=&  \sum_{k \geq 0} \frac{1}{2-p} h_t(k-1)^{2-p} g_t(k-1)^{2p-2}+ \frac{1-p}{2-p} g_t(k)^2 f_t(k)^{p-2}\\
&& -2h_t(k-1) f_t(k)^{p-1} \\ 
&&+ \frac{1}{2-p} h_t(k-1)^{2-p} g_t(k)^{2p-2} + \frac{1-p}{2-p} g_t(k-1)^2 f_t(k)^{p-2} .
\end{eqnarray*}

We denote $v_{+,t}(k):= \frac{g_t(k)}{f_t(k)}$ and $v_{-,t}(k) :=\frac{g_t(k-1)}{f_t(k)}$. The Benamou-Brenier equation~\eqref{eq:BBCondZ} is then written $h_t(k-1)=v_{+,t}(k)v_{-,t}(k) f_t(k)$ and we have

\begin{eqnarray*}
\sum_{k \geq 0} \nabla_2(h_t(k)) f_t(k)^{p-1} &\leq& \sum_{k \geq 0} \frac{1}{2-p} v_{+,t}(k)^{2-p} v_{-,t}(k)^p f_t(k)^p + \frac{1-p}{2-p} v_{+,t}(k)^2 f_t(k)^p\\
&& - 2 v_{+,t}(k)v_{-,t}(k) f_t(k)^p + \frac{1}{2-p} v_{+,t}(k)^p v_{-,t}(k)^{2-p} f_t(k)^p + \frac{1-p}{2-p} v_{-,t}(k)^2 f_t(k)^p.
\end{eqnarray*}

With the same notations, we have
\begin{equation}
(g_t(k)-g_t(k-1))^2 f_t(k)^{p-2} = (v_{+,t}(k)^2-2v_{+,t}(k)v_{-,t}(k)+v_{-,t}(k)^2) f_t(k)^{p},
\end{equation}
and
\begin{equation}
0 = \sum_{k \geq 0} g_t(k)^{p}-g_t(k-1)^{p} = \sum_{k \geq 0} (v_{+,t}(k)^p-v_{-,t}(k)^p)f_t(k)^p.
\end{equation}

We use these estimations and the positivity of $\psi$ to write
\begin{eqnarray*}
-\frac{1}{p} H_p''(t) &=& \sum_{k \geq 0} \nabla_2(h_t(k)) f_t(k)^{p-1}-(1-p) \sum_{k\geq 0} \nabla (g_t(k))^2 f_t(k)^p \\
&& -\frac{1}{2-p} \sum_{k \geq 0} g_t(k)^{p}-g_t(k-1)^{p} \\
&\leq& \sum_{k \geq 0} v_{+,t}(k)^2 f_t(k)^p \psi\left(\frac{v_{-,t}(k)}{v_{+,t}(k)}\right) \\
&\leq& 0,
\end{eqnarray*}
which finishes the proof of the theorem. $\square$

\subsection{Binomial mixtures and $W_2$-optimal couplings.}

In this appendix we consider two finitely supported probability measures $f_0$ and $f_1$ on $\mathbb{Z}$. Through this article and the previous one (see~\cite{HillionGeodesic}), we have seen that, by considering a mixture of binomial measures with respect to a proper coupling between $f_0$ and $f_1$, it is possible to construct a $W_1$-geodesic $(f_t)_{t \in [0,1]}$ which satisfies a Benamou-Brenier condition~\eqref{eq:BBcondition} which is a discrete analogue of a characterization of $W_2$-geodesics on the real line. 

\medskip

Another natural way to generalize the notion of $W_2$-geodesic from the continuous setting to the discrete setting is the following:

\begin{defi}
Let $\pi$ be the unique $W_2$-optimal coupling between $f_0$ and $f_1$. The binomial/$W_2$ interpolation $(f_t)_{t \in [0,1]}$ is defined by:
\begin{equation} 
f_t(k) : =\sum_{(i,j)} \pi_{i,j} \bino_{(i,j),t}(k), 
\end{equation}
where $\bino_{(i,j),t}$ is the binomial family between $i$ and $j$.
\end{defi}

Basic theorems on optimal transportation give the existence and uniqueness of a $W_2$-optimal coupling $\pi$ between $f_0$ and $f_1$. Thus the binomial/$W_2$ interpolation $(f_t)_{t \in [0,1]}$ exists and is unique.

\medskip

The question of the convexity of the entropy along $(f_t)$ is still open. The particular case where $f_1$ is a translation of $f_0$ has been studied by the author in~\cite{HillionTranslation}. In this appendix we prove the more general:

\begin{thm}\label{th:EntroBinoW2} We make the following assumptions:
\begin{enumerate}
 \item The measure $f_0$ is stochastically dominated by $f_1$ : $f_0 << f_1$, which means that for each $l \in \mathbb{Z}$, $\sum_{l \leq k} f_0(l) \geq \sum_{l \leq k} f_1(l)$.
 \item Each $f_t$ is log-concave, i.e. that the inequality $f_t(k+1)^2 \geq f_t(k)f_t(k+2)$ holds for any $t \in [0,1]$ and $k \in \mathbb{Z}$.
\end{enumerate}
Then the entropy $H(t)$ of $f_t$ is a convex function of $t$.
\end{thm}

The stochastic domination assumption is not necessary but allows us to give a simpler proof. We will use it through the following:
\begin{lem}\label{lem:W2coupling}
We suppose that $f_0 << f_1$. Then the $W_2$-optimal coupling $\pi$ between $f_0$ and $f_1$ satisfies the following:
\begin{itemize}
 \item If $\pi(i,j)>0$ then $i \leq j$.
 \item If $\pi(i_1,j_1)>0$ and $\pi(i_2,j_2)>0$ then $(i_2-i_1)(j_2-j_1) \geq 0$.
\end{itemize}
\end{lem}

\textbf{Remark.} In particular the stochastic domination assumption allows us to use the same notations $g(k):=g(k,k+1)$, $h(k) := h(k,k+1,k+2)$ as in the first part of the Appendix.

\medskip

\textbf{Proof of Theorem~\ref{th:EntroBinoW2}: } Using the first point of Lemma~\ref{lem:W2coupling}, we can write
\begin{equation}
f_t(k) = \sum_{i \leq j} \pi(i,j) \bino_{(j-i),t}(k-i). 
\end{equation}
We now define the families of functions $(g_t)_{t \in [0,1]}$ and $(h_t)_{t \in [0,1]}$ by:
\begin{eqnarray*}
g_t(k) &:=&  \sum_{i \leq j} \pi(i,j) (j-i) \bino_{(j-i),t}(k-i) \\ h_t(k) &:= & \sum_{i \leq j} (j-i)(j-i-1) \pi(i,j) \bino_{(j-i-2),t}(k-i),
\end{eqnarray*}
so we have
\begin{equation}
\frac{\partial}{\partial t} f_t(k) = -\nabla g_t(k) \ , \  \frac{\partial}{\partial t} g_t(k) = -\nabla h_t(k).
\end{equation}
The study of the entropy of $f_t$ is similar to case of $W_{1,+}$-interpolations. We have:
\begin{eqnarray*}
H''(t) &=& \sum_{k \in \mathbb{Z}} \nabla_2h_t(k) \log(f_t(k)) + \sum_{k \in \mathbb{Z}} \frac{(\nabla g_t(k))^2}{f_t(k)}.
\end{eqnarray*}
The major difference with $W_{1,+}$-interpolations comes from the fact that $f_t(k)h_t(k-1)$ is a priori not equal to $g_t(k)g_t(k-1)$. Let us introduce the family of functions $\tilde{h_t}(k) := \frac{g_t(k)g_t(k+1)}{f_t(k+1)}$. The triple $(f_t,g_t,\tilde{h}_t)$ is a BB-triple on $\mathbb{Z}$, so we have:
\begin{eqnarray*}
H''(t) &=& \sum_{k \in \mathbb{Z}} \nabla_2(h_t-\tilde{h}_t)(k) \log(f_t(k))+\sum_{k \in \mathbb{Z}} \nabla_2\tilde{h}_t(k) \log(f_t(k)) + \sum_{k \in \mathbb{Z}} \frac{(\nabla g_t(k))^2}{f_t(k)} \\
&\geq& \sum_{k \in \mathbb{Z}} \nabla_2(h_t-\tilde{h}_t)(k) \log(f_t(k)) \\
&=& \sum_{k \in \mathbb{Z}} (h_t-\tilde{h}_t)(k) \nabla_2 \log(f_t(k+2)).
\end{eqnarray*}
By the assumption on the log-concavity of $f_t$, it thus suffices to show that $h_t \leq \tilde{h}_t$. to prove this fact, we notice that we can write $g_t(k)$, $g_t(k-1)$ and $h_t(k-1)$ under the form:
\begin{equation*}
g_t(k) = \sum_{i \leq j} \pi(i,j) \bino_{(j-i),t}(k-i) \frac{j-k}{1-t}  \ , \ g_t(k) = \sum_{i \leq j} \pi(i,j) \bino_{(j-i),t}(k-i) \frac{k-i}{t}, 
\end{equation*}
\begin{equation*}
h_t(k-1) = \sum_{i \leq j} \pi(i,j) \bino_{(j-i),t}(k-i) \frac{(j-k)(k-i)}{t(1-t)}.
\end{equation*}

Let us denote, for $i \leq j$, $a(i,j) := \pi(i,j) \bino_{(j-i),t}(k-i)$. Then $g_t(k)g_t(k-1)-f_t(k)h_t(k-1)$  can be seen as a quadratic form in the variables $(a(i,j))_{(i,j) \in \supp(\pi)}$. The coefficient associated to $a(i,j)^2$ is
\begin{equation}
\frac{j-k}{1-t} \frac{k-i}{t} - \frac{(j-k)(k-i)}{t(1-t)} =0 .
\end{equation}
If $(i_1,j_1) \neq (i_2,j_2)$ are in $\supp(\pi)$ then the coefficient associated to $a(i_1,j_1) a(i_2,j_2)$ is
\begin{eqnarray*}
\frac{j_1-k}{1-t}\frac{k-i_2}{t}+\frac{j_2-k}{1-t}\frac{k-i_1}{t}-\frac{(j_1-k)(k-i_1)}{t(1-t)}-\frac{(j_2-k)(k-i_2)}{t(1-t)}
&=& \frac{(j_2-j_1)(i_2-i_1)}{t(1-t)} \\&\geq& 0.
\end{eqnarray*}

This shows that $h_t \leq \tilde{h}_t$, and finishes the proof of Theorem~\ref{th:EntroBinoW2}. $\square$

\end{document}